\documentclass[12pt]{article}

\usepackage{graphicx, amsmath}
\usepackage{amsfonts}
\usepackage{verbatim}
\usepackage{latexsym}
\usepackage{cite}

\usepackage{amssymb}
\usepackage{color,soul}
\usepackage{booktabs}
\usepackage{multirow}
\usepackage{tabularx}
\usepackage{multirow}
\usepackage{graphicx}
\usepackage{bicaption}
\usepackage{amsmath}
\usepackage{ntheorem}

\makeatletter
\@addtoreset{equation}{section}

\numberwithin{equation}{section}
\begin{document}
	
	\newcommand{\E}{\mathbb{E}}
	\newcommand{\PP}{\mathbb{P}}
	\newcommand{\RR}{\mathbb{R}}

	\newtheorem{theorem}{Theorem}[section]
	\newtheorem{lemma}{Lemma}[section]
	\newtheorem{coro}{Corollary}[section]
	\newtheorem{defn}{Definition}[section]
	\newtheorem{assp}{A}
	\newtheorem{expl}{Example}[section]
	\newtheorem{prop}{Proposition}[section]
	\newtheorem{remark}{Remark}[section]
	\theoremseparator{:}
	\newtheorem{proof}{Proof}
	\newtheorem*{Proof}{Proof} 
	\newtheorem{PROOF}{Proof of Theorem}[section]
	\newcommand\tq{{\scriptstyle{3\over 4 }\scriptstyle}}
	\newcommand\qua{{\scriptstyle{1\over 4 }\scriptstyle}}
	\newcommand\hf{{\textstyle{1\over 2 }\displaystyle}}
	\newcommand\hhf{{\scriptstyle{1\over 2 }\scriptstyle}}
	
	\newcommand{\eproof}{\indent\vrule height6pt width4pt depth1pt\hfil\par\medbreak}
	
	\newcommand\wD{\widehat{\D}}
	\newcommand{\ka}{\kappa_{10}}
	
	\title{
		\bf {Convergence of Numerical Solution of The Tamed Milstein Method for NSDDEs}
	}
	\author
{ {\bf Qiquan Fang}$^{\tt a}$\ \, {\bf Yingxiao Min }$^{\tt a}$ \ \,{\bf Yingying Wang }$^{\tt a}$\ \, {\bf Yanting Ji }$^{\tt a}\thanks{Contact e-mail address: yanting.ji@zust.edu.cn}$
	\\ [1ex]$^{\tt a}$  Department of Mathematics\\
	Zhejiang University of Science and Technology, Hangzhou, China
	\\}

	\date{}
	
	\maketitle
	\begin{abstract}
	In this paper, we apply the tamed technique to the Milstein numerical scheme to investigate Neutral Stochastic Delay Differential Equations(NSDDEs) with highly nonlinear coefficients. Under the local Lipschitz condition and Khasminskii condition, the tamed Milstein numerical solution converges strongly to the exact solution.
		\\

		\smallskip
		
		{\it Key words}: neutral stochastic differential delay equations; tamed Milstein scheme; local Lipschitz; strong convergence.
		
	\end{abstract}
	\section{Introduction}

As an important class of differential dynamical systems, Neutral Stochastic Delay  Differential Equations(NSDDEs) play an important role in many fields, such as biology, finance and automatic control\cite{2014State,Guo2018The}. In general, although some linear stochastic differential equations(SDEs) can be solved explicitly, solutions to SDEs with highly nonlinear coefficients are difficult to obtain. Therefore, numerical schemes are of great importance in the study of nonlinear SDEs. In recent years, many literatures related to numerical schemes have emerged, see\cite{1955Continuous,2008Mean}.

Most convergence results on the numerical schemes for SDEs require the coefficients satisfying the local Lipschitz condition and linear growth condition.
When the drift term or diffusion term does not meet the linear growth condition, the Euler-Maruyama(EM) method will divergence in finite time, see\cite{2011Strong}. Later, Hutzenthaler\cite{2012STRONG} proposed a tamed EM method with strong convergence of order $1/2$. Then, Sabanis\cite{2013A} extended the strong convergence results of the tamed Euler scheme to SDEs whose drift term satisfies local conditions and superlinear growth conditions. Wang and Gan\cite{2013The} further extended a tamed Milstein scheme to SDEs with commutative noise.

NSDDEs is one of the special class of ordinary SDEs, which can be used to describe important phenomena in life. When their coefficients do not meet the linear growth, the classical  EM scheme does not converge. Therefore, scholars have given many modified numerical schemes for NSDDEs with highly nonlinear coefficients, and studied their convergence. Ji\cite{2017Tamed} analyzed the convergence rate of tamed EM for NSDDEs with linear growth of diffusion term. Tan and Yuan\cite{2019Convergence} studied the strong convergence of the tamed theta scheme for NSDDEs with one-side Lipschitz drift coefficients. Yan et al.\cite{2017Strong} investigated the strong convergence of the split-step theta method for the NSDDEs, where the corresponding coefficients may be highly nonlinear with respect to the delay variables.  Lan and Wang\cite{2019Strong} investigated strong convergence of modified truncated EM method for NSDDEs. For a class of highly nonlinear SDDEs with nonlinear growth conditions, Zhang et al.\cite{2019Zhang} proposed the truncated Milstein method and gave its convergence rate. 
Deng et al. \cite{2021Tamed} developed two types of explicit tamed EM schemes for NSDDEs, in which both drift and diffusion coefficients can grow superlinearly, and investigate the strong convergence, mean-square stability.

Inspired by \cite{2013The,2017Tamed,2021Tamed}, we develop tamed Milstein scheme for NSDDEs with highly nonlinear coeffcients. The rest of the paper is organized as follows. In the section $\ref{preliminaries}$, we present some preliminaries and assumptions on the NSDDEs. The tamed Milstein scheme is proposed in section $\ref{tamed milstein shceme}$ and the strong convergence of the tamed Milstein is given in section $\ref{convergence of the tamed milstein}$.
\section{Preliminaries}\label{preliminaries}
Throughout this paper, unless otherwise specified, let $(\Omega,\mathcal{F},P)$ be a complete probability space with a filtration $\left\{\mathcal{F}_{t}\right\}_{t \geq 0}$ satisfying the usual conditions(i.e., it is right continuous and increasing while $\mathcal{F}_0$ contains all P-null sets). Let $B\left(t\right)$ be an m-dimensional Brownian motion. $\|\cdot\|$denotes the Euclidean norm in $\mathbb{R}^n$. The inner product of $x, y$ in $\mathbb{R}^n$ is denoted by $\langle x, y\rangle$ or $x^{T} y$. If A is a matrix, its trace norm  is denoted by $|A|=\sqrt{\operatorname{trace}\left(A^{T} A\right)}$. Let $\tau>0$ be a constant and denote $C\left(\left[-\tau,0\right];\mathbb{R}^n\right)$ the space of all continuous functions from $\left[-\tau,0\right]$ to $\mathbb{R}^n$ with the norm $\|\xi\|=\sup_{-\tau \leq \theta \leq 0}|\xi(\theta)|$.

In this paper, we study the numerical approximation of the neutral stochastic differential delay equations(NSDDEs)
\begin{equation}\label{NSDDEs1}
d[x(t)-D(x(t-\tau))]=b(x(t), x(t-\tau)) d t+\sigma(x(t), x(t-\tau)) d B(t), t \geq 0,
\end{equation}
with the initial data 
\begin{equation}\label{initial data}
	\{x(\theta):-\tau \leq \theta \leq 0\}=\xi \in L_{\mathcal{F}_{0}}^{p}\left([-\tau, 0] ; \mathbb{R}^{n}\right),
\end{equation}
that is $\xi$ is an $\mathcal{F}_{0}$-measurable $C\left(\left[-\tau,0\right];\mathbb{R}^n\right)$-valued random variable and $\mathbb{E}\|\xi\|^p<\infty$.
Here
\begin{equation*}
	D:\mathbb{R}^{n} \rightarrow \mathbb{R}^{n} ,\ b: \mathbb{R}^{n} \times \mathbb{R}^{n} \rightarrow \mathbb{R}^{n} ,\    \sigma: \mathbb{R}^{n} \times \mathbb{R}^{n} \rightarrow \mathbb{R}^{n \times m}.
\end{equation*}

We assume that the coefficients $D,\  b$ and $\sigma$ are Borel-measurable and satisfy the following conditions: \\
\textbf{(A1)} For any $s,t\in[-\tau,0]$ and $q>0$, there exists a positive constant $\theta>0$ such that
\begin{equation*}
	\mathbb{E}\left|\xi\left(t\right)-\xi\left(s\right)\right|^q\le\theta\left|t-s\right|^q.
\end{equation*}
\textbf{(A2)}  $D(0)=0$ and there exists a constant $\kappa \in(0,1)$ such that
\begin{equation*}
	|D\left(x\right)-D\left(y\right)\left|\le \kappa \right|x-y| \ \rm for\ all \ \emph x,\emph y\in \mathbb{R}^n.
\end{equation*}
\textbf{(A3)} For every constant $R>0$, there exists a positive constant $K_{R}$,  such that
\begin{equation*}
|b(x, y)-b(\bar{x}, \bar{y})| \vee|\sigma(x, y)-\sigma(\bar{x}, \bar{y})| \leq K_{R}\left(|x-\bar{x}|+|y-\bar{y}|\right),
\end{equation*}
for all $|x| \vee|y| \vee|\bar{x}| \vee|\bar{y}| \leq R$.\\
\textbf{(A4)} There are constants $R>0$ and $\bar{K}_R>0$ such that
\begin{align*}
	|\sigma_i(x,y)\sigma(x,y)-\sigma_i(\bar{x},\bar{y})\sigma(\bar{x},\bar{y})|\leq \bar{K}_{R}(|x-\bar{x}|+|y-\bar{y}|),\quad i=1,2.
\end{align*}
for all $|x| \vee|y| \vee|\bar{x}| \vee|\bar{y}| \leq R$ where 
\begin{align*}
 \sigma_{i}\left(x_{1}, x_{2}\right)=\frac{\partial \sigma\left(x_{1}, x_{2}\right)}{\partial x_{i}}.
\end{align*} 
\textbf{(A5)} There exist constants $p>$2 and $K_{1}>0$ such that 
\begin{equation*}
	\left(x-D(y)\right)^Tb(x, y)+\frac{p-1}{2}\left|\sigma(x, y)\right|^2\le K_{1}(1+\left|x\right|^2+{|y|}^2),
\end{equation*}
for  $x, y\in \mathbb{R}^n$.
\begin{lemma}\label{x(t)}
	Let assumption \rm\textbf{(A3)} \emph{and} \rm\textbf{(A4)} \emph{hold, the equation\eqref{NSDDEs1} have a unique global solution $x(t),t \in [0,T]$, and there exists a positive constant $C$ such that for any $p\geq2$ }
	\begin{equation}
		\mathbb{E}\Big(\sup _{0\leq t \leq T}\big|x(t)\big|^p)\Big)\leq C,
			\end{equation}
			where the positive constant $C:=C(p,T,K_1,||\xi||,\kappa).$
\end{lemma}

	\begin{Proof}
		 \rm By  ${(a+b)}^p\le\left(1+\varepsilon\right)^{p-1}\left(a^p+\varepsilon^{1-p}b^p\right)$ where $a,b,\varepsilon>0,p\geq2$
\begin{align*}
	|x(t)|^{p}&=|x(t)-D(x(t-\tau))+D(x(t-\tau))|^{p} \\
	&\leq(1+\varepsilon)^{p-1}\Big(|x(t)-D(x(t-\tau))|^p+\varepsilon^{1-p}|D(x(t-\tau))|^p\Big).
\end{align*}
	Let $\varepsilon=\frac{\kappa}{1-\kappa}$, then
	\begin{align}\label{basic inequality}
		|x(t)|^p&\leq(1-\kappa)^{1-p}|x(t)-D(x(t-\tau))|^{p}+\kappa^{1-p}D(x(t-\tau))^p\nonumber\\
		&\leq(1-\kappa)^{1-p}|x(t)-D(x(t-\tau))|^{p}+\kappa x(t-\tau)^p\nonumber\\
		&\leq(1-\kappa)^{1-p}|x(t)-D(x(t-\tau))|^{p}+\kappa(||\xi||^p+x(t))^p)\nonumber\\
		&\leq(1-\kappa)^{-p}|x(t)-D(x(t-\tau))|^p+\frac{\kappa}{1-\kappa}||\xi||^p.
	\end{align}
An application of $
\text { Itô }
$ formula yields
	\begin{align}\label{Ito}
		&|x(t)-D(x(t-\tau))|^{p}\nonumber\\
		&\leq\big|\xi(0)-D(\xi(-\tau))\big|^p+\int_{0}^{t}p\big|x(s)-D(x(s-\tau))\big|^{p-2}\Big[\frac{p-1}{2}\big|\sigma(x(s),x(s-\tau)|^2\nonumber\\
		&\quad+\big(x(s)-D(x(s-\tau))\big)^T b\big(x(s),x(s-\tau)\big)\Big]ds
		\nonumber\\
		&+\int_{0}^{t}p\big|x(s)-D(x(s-\tau))\big|^{p-2}\big(x(s)-D(x(s-\tau))\big)^T \sigma\big(x(s),x(s-\tau)\big)dB(s) \nonumber\\
		&=\big|\xi(0)-D(\xi(-\tau))\big|^p+I_{1}(t)+I_{2}(t)+I_{3}(t).
	\end{align}
	Under \textbf{(A2)} and \textbf{(A5)}, we have the following estimate  
	\begin{align}\label{sum of I1 and I2}
		\mathbb{E}&\big(\sup _{0 \leq t\leq T}I_{1}(t)\big)+\mathbb{E}\big(\sup _{0 \leq t\leq T}I_{2}(t)\big)\nonumber\\
		=&\mathbb{E}\Big(\sup_{0\leq t\leq T}\int_{0}^{t}p\big|x(s)-D(x(s-\tau))\big|^{p-2}\Big[\big(x(s)-D(x(s-\tau))\big)^T b\big(x(s),x(s-\tau)\big)
		\nonumber\\
		&+\frac{p-1}{2}\big|\sigma(x(s),x(s-\tau)\big|^2\Big]ds\Big)\nonumber\\
		\leq &C\mathbb{E}\int_{0}^{T}\big(x(s))^{p-2}+|D(x(s-\tau))|^{p-2})\left(1+|x(s)|^{2}+|x(s-\tau)|^{2}\right)ds\nonumber \\
		\leq &C\mathbb{E}\int_{0}^{T}\Big(1+\big|x(s)\big|^p+\big|x(s-\tau)\big|^p\Big)ds\nonumber\\
	\leq & C+C\int_{0}^{T}\mathbb{E}\Big(\sup _{0\leq u\leq s}\big|x(u)\big|^p\Big)ds.
	\end{align}
	For $I_{3}(t)$, by the Burkholder-Davis-Gundy(BDG) inequality and H$\ddot{\boldsymbol{\text{o}}}$lder inequality which yields                                                                                                                                                                                                                                                                                                                                                                                                                                                                                                                                                                                                                                                                                                                            
	\begin{align*}\label{I3}
	\mathbb{E}&\big(\sup _{0\leq t\leq T}I_{3}(t)\big)\nonumber\\
	=&\mathbb{E}\Bigg(\sup_{0\leq t\leq T}\int_{0}^{t}p\big|x(s)-D(x(s-\tau))\big|^{p-2}\big(x(s)-D(x(s-\tau))\big)^T \sigma\big(x(s),x(s-\tau)\big)dB(s) \Bigg)\nonumber\\
	\leq &C_{p}\mathbb{E}\Bigg(\int _{0}^{T}\Big(\big|x(s)-D(x(s-\tau))\big|^{2p-2}\big|\sigma(x(s),x(s-\tau))\big|^2\Big)ds \Bigg)^\frac{1}{2}.\nonumber \\
	\leq &C_{p}\mathbb{E}\Bigg(\sup _{0\leq t\leq T}\Big|x(t)-D\big(x(t-\tau)\big)\Big|^{p-1}\Big(\int _{0}^{T}\big|\sigma(x(s),x(s-\tau))\big|^2ds\Big)^{1/2}\Bigg).
	\end{align*}
	Then by inequality $a^{p-1}b\leq\frac{p-1}{p}a^{p}+\frac{1}{p}b^p$ and \textbf{(A3)}, we could obtain
	\begin{align}
	\mathbb{E}&\big(\sup _{0\leq t\leq T}I_{3}(t)\big)\nonumber\\
	\leq &C_{p} \mathbb{E}\Big(\sup _{0\leq t\leq T}\Big|x(t)-D\big(x(t-\tau)\big)\Big|^{p}\Big)+C_{p}\mathbb{E}\Big(\int _{0}^{T}|\sigma(x(s),x(s-\tau))\big|^2ds\Big)^{p/2}\nonumber \\
	\leq &C_{p} \mathbb{E}\Big(\sup _{0\leq t\leq T}\Big|x(t)-D\big(x(t-\tau)\big)\Big|^{p}\Big)+C_{p}\mathbb{E}\Big(\int _{0}^{T}(1+|x(s)|^2+|x(s-\tau)|^2ds\Big)^{p/2}\nonumber \\
	\leq &C_{p}+C_{p} \mathbb{E}\Big(\sup _{0\leq t\leq T}\Big|x(t)-D\big(x(t-\tau)\big)\Big|^{p}\Big)+C_{p}\int_{0}^{T}\mathbb{E}\Big(\sup _{0\leq u\leq s}\big|x(u)\big|^p\Big)ds.
	\end{align}
	Now, substituting \eqref{sum of I1 and I2} and \eqref{I3} into \eqref{Ito}, it derives that
	\begin{equation}
		\mathbb{E}\Big(\sup _{0\leq t \leq T}\big|x(t)-D(x(t-\tau))\big|^p\Big)\leq C+C\int_{0}^{T}\mathbb{E}\Big(\sup _{0\leq u\leq s}\big|x(u)\big|^p\Big)ds.
	\end{equation}
	Then, substituting into \eqref{basic inequality}, it derives that
	\begin{equation*}
		\mathbb{E}\big(\sup _{0\leq t\leq T}|x(t)|^{p}\big)\leq C+C\int_{0}^{T}\mathbb{E}\Big(\sup _{0\leq u\leq s}\big|x(u)\big|^p\Big)ds.
	\end{equation*}
	An application of Gronwall inequality yields
	\begin{equation*}
		\mathbb{E}\Big(\sup _{0\leq t \leq T}\big|x(t)\big|^p)\Big)\leq C.
			\end{equation*}
	\end{Proof}

\section{Tamed Milstein Scheme}\label{tamed milstein shceme}
Fix $0<\tau <T$, without loss of gengrality, we assume that $T$ and $\tau$ are rational numbers, and the step size $\Delta \in(0,1)$ be fraction of $\tau$ and $T$, so that there exist two positive integers $M$ and $m$ such that $\Delta=T/M=\tau/m,t_{k}=k \Delta: k=-m, \ldots 0,1,2, \ldots M-1$.

Applying the classic Milstein scheme\cite{G1974Approximate} to \eqref{NSDDEs1} with initial data 
\eqref{initial data}, we have the discrete-time Milstein shceme
\begin{equation}
\left\{\begin{aligned}
Y_{t_{k}}=& \xi\left(t_{k}\right), \quad k=-m,-m+1, \ldots, 0 \\
Y_{t_{k+1}}=&D(Y_{t_{k+1-m}})+Y_{t_{k}}-D\left(Y_{t_{k-m}}\right)+b_{h}\left(Y_{t_{k}}, Y_{t_{k-m}}\right) \Delta
+\sigma\left(Y_{t_{k}}, Y_{t_{k-m}}\right)
 \Delta B_{t_{k}}\\
 &+\sigma_{1}\left(Y_{t_{k}}, Y_{t_{k-m}}\right) \sigma\left(Y_{t_{k}}, Y_{t_{k-m}}\right) l_{1}\\
 &+\sigma_{2}\left(Y_{t_{k}}, Y_{t_{k-m}}\right) \sigma\left(Y_{t_{k-m}}, Y_{t_{k-2m}}\right) l_{2},\quad k=0,1, \ldots, M-1,
\end{aligned}\right.
\end{equation}
Here
\begin{align*}
	l_{1}=\int_{t_{k}}^{t_{k+1}} \int_{t_{k}}^{s} d B(t) d B(s)=\frac{\left(\Delta B_{t_{k}}\right)^{2}-\Delta}{2}, \ l
_{2}=\int_{t_{k}}^{t_{k+1}} \int_{t_{k}}^{s} d B(t-\tau) d B(s).
\end{align*}
\begin{align*}
	\Delta B_{t_{k}}=B_{t_{k+1}}-B_{t_{k}},  \ \sigma_1\left(x,y\right)=\frac{\partial\sigma\left(x,y\right)}{\partial x}, \ \sigma_2\left(x,y\right)=\frac{\partial\sigma\left(x,y\right)}{\partial y}.\end{align*}
Defining the tamed drift term
\begin{align}
	b_{h}(x,y)=\frac{b(x,y)}{1+\Delta^{\alpha}\big|b(x,y)\big|},
\end{align}
for all $x,y\in \mathbb{R}^n$ and $\alpha \in\left(0, \frac{1}{2}\right] .$
\\By observation, one has
\begin{align}\label{bh}
	|b_h(x,y)|\leq {\frac{1}{\Delta^{\alpha}} \wedge |b(x,y)|}.
\end{align}
\begin{remark}\label{remark}
	When $b$ is replaced by $b_h$, assumptions \rm{\textbf{A3}} \emph {and} \rm{\textbf{A4}} \emph{are still true}\\
Under assumption \textbf{(A3)}
\begin{align}\label{bh1}
	&|b_{h}(x, y)-b_{h}(\bar{x}, \bar{y})|\nonumber\\
	&=\Big|\frac{b(x,y)}{1+\Delta^{\alpha}|b(x,y)|}-\frac{b(\bar{x},\bar{y})}{1+\Delta^{\alpha}|b(\bar{x},\bar{y})|} \Big|\nonumber\\
	&=\frac{b(x,y)-b(\bar{x},\bar{y})}{(1+\Delta^{\alpha}|b(x,y)|)(1+\Delta^{\alpha}|b(\bar{x},\bar{y})|)}+\frac{\Delta^{\alpha}[b(x,y)|b(\bar{x},\bar{y})|-|b(x,y)|b(\bar{x},\bar{y})]}{(1+\Delta^{\alpha}|b(x,y)|)(1+\Delta^{\alpha}|b(\bar{x},\bar{y})|)}\nonumber\\
	&\leq\frac{1}{(1+\Delta^{\alpha}|b(x,y)|)(1+\Delta^{\alpha}|b(\bar{x},\bar{y})|)}K_{R}\left(|x-\bar{x}|+|y-\bar{y}|\right)\nonumber\\
	&\le K_{R}\left(\left|x-\bar{x}\right|+\left|y-\bar{y}\right|\right).
\end{align}
Under assumption \textbf{(A4)}
\begin{align}\label{yuan B2}
	&\left(x-D(y)\right)^Tb_h(x,y)\nonumber\\
	&=\frac{1}{1+\Delta^{\alpha}|b(x,y)|}\left(x-D(y)\right)^Tb(x,y)\nonumber\\
	&\leq \frac{K_{1}}{1+\Delta^{\alpha}|b(x,y)|}(1+|x|^2+|y|^2)\nonumber\\
	&\leq K_{1}(1+|x|^2+|y|^2).
\end{align}
\end{remark}
Define tamed Milstein scheme continuous-time step process as 
\begin{equation}
\bar{y}(t)=\sum_{k=0}^{\infty} Y_{t_{k}} I_{\left[t_{k}, t_{k+1}\right)}(t), t \geq 0.
\end{equation}
where $I_{\left[t_{k}, t_{k+1}\right)}(t)$is the indicator function on $[t_k,t_{k+1})$.\\
Then, we could define the corresponding continuous-time tamed Milstein scheme
\begin{equation}
\left\{\begin{aligned}\label{tamed method}
y\left(t\right)&=\xi(t), t\in\left[-\tau,0\right];\\
y(t)&=D(\bar{y}(t-\tau))+\xi(0)-D(\xi(-\tau))\\
&\quad+\int_{0}^{t} b_{h}(\bar{y}(s), \bar{y}(s-\tau)) d s+\int_{0}^{t} \sigma(\bar{y}(s), \bar{y}(s-\tau)) d B(s) \\
&\quad+\int_{0}^{t} \sigma_{1}(\bar{y}(s), \bar{y}(s-\tau)) \sigma(\bar{y}(s), \bar{y}(s-\tau)) \Delta \bar{B}(s) d B(s)  \\
& \quad+\int_{0}^{t} \sigma_{2}(\bar{y}(s), \bar{y}(s-\tau)) \sigma(\bar{y}(s-\tau), \bar{y}(s- 2\tau)) \Delta \bar{B}(s-\tau) d B(s),t\in(0,T].
\end{aligned}\right.
\end{equation}
Here$\ \Delta \bar{B}(t)=\sum_{k=0}^{\infty}\left(B_{t_{k+1}}-B_{t_{k}}\right) I_{\left[t_{k}, t_{k+1}\right)}(t), t \geq 0	$.

Then $\bar{y}(t)=Y_{t_{k}}=y(t_k)$, for any $t\in[t_k,t_{k+1})$ with $k\geq 0$, we could conclude that
\begin{align}
	\mathbb{E}|\bar{y}(t)|\leq \mathbb{E}\big(\sup_{0\leq t\leq T}|y(t)|\big).
\end{align}
\textbf{(B1)}
There exists a positive constant $N_R$ depending on R, such that for any $x,y\in \mathbb{R}^d,$
\begin{align}
	\sup_{|x|\vee|y|\leq R}|b(x,y)-b_h(x,y)|\leq N_R\Delta.
\end{align}

\section{Convergence of the tamed Milstein }\label{convergence of the tamed milstein}
In this section, we show that the tamed Milstein scheme converges to the exact solution under certain conditions, i.e. we have the following main result:
\begin{theorem}
	Let \rm\textbf{(A1)} \textbf{(A2)} \emph{remark 3.1 and \rm\textbf{(B1)} \emph {hold, then for any $p\geq2$}
	\begin{equation}
\lim _{\Delta \rightarrow 0} \mathbb{E}\left[\sup _{0 \leq t \leq T}|x(t)-y(t)|^{p}\right]=0.
\end{equation}}
\end{theorem}
\subsection{Moment properties of tamed Milstein method}
To prove our main results, in this subsection we investigate the boundedness of moments tamed Milstein approximation in the following section.
\begin{lemma}
	Under \eqref{yuan B2}, for$\ all\ t\in[0,\ T]$, we have
	\begin{equation}\label{y(t)}
\mathbb{E}\sup _{0<t \leq T} |y(t)|^{p} \leq C,
\end{equation}
where the positive constant $C:=C(p,T,K_R,K_1,||\xi||,\kappa).$
\end{lemma}
\begin{Proof}
	\rm  Using the similar way as in \eqref{basic inequality},
	\begin{equation}\label{yt}
	\left|y\left(t\right)\right|^p\le(1-\kappa)^{-p}{|y\left(t\right)-D\left(\bar{y}\left(t-\tau\right)\right)|}^p+\frac{\kappa}{1-\kappa}{|D\left(\bar{y}\left(t-\tau\right)\right)|}^p.
\end{equation}
By applying the $
\text { Itô }
$ formula to \eqref{tamed method}, it derives that 
\begin{align}\label{sum of H}
	&|y\left(t\right)-D\left(\bar{y}\left(t-\tau\right)\right)|^p\nonumber\\
	&\leq \left|\xi(0)-D(\xi(-\tau))\right|^{p}\nonumber \\
	&+\int_{0}^{t} p|y(s)-D(\bar{y}(s-\tau))|^{p-2}\left[\left(\bar{y}(s)-D(\bar{y}(s-\tau))^{T}  b_{h}(\bar{y}(s), \bar{y}(s-\tau))\right.\right.\nonumber  \\
	&\quad+(p-1)\big|\sigma(\bar{y}(s), \bar{y}(s-\tau))\big|^2+(p-1)\big|\sigma_{1}(\bar{y}(s), \bar{y}(s-\tau)) \sigma(\bar{y}(s), \bar{y}(s-\tau)) \Delta \bar{B}(s)\nonumber \\
	&\quad+\sigma_{2}(\bar{y}(s), \bar{y}(s-\tau)) \sigma(\bar{y}(s), \bar{y}(s-\tau)) \Delta \bar{B}(s-\tau)\big|^2\nonumber \\
	&\quad\left.+(y(s)-\bar{y}(s))^{T} b_{h}(\bar{y}(s), \bar{y}(s-\tau))\right] ds\nonumber \\
	&+\int_{0}^{t} p|y(s)-D(\bar{y}(s-\tau))|^{p-2}\big(y(s)-D(\bar{y}(s-\tau))\big)^T \Big| \sigma(\bar{y}(s), \bar{y}(s-\tau))\nonumber \\
	&\quad+\sigma_{1}(\bar{y}(s), \bar{y}(s-\tau)) \sigma(\bar{y}(s), \bar{y}(s-\tau)) \Delta \bar{B}(s)\nonumber \\
	&\quad+\sigma_{2}(\bar{y}(s), \bar{y}(s-\tau)) \sigma(\bar{y}(s-\tau), \bar{y}(s-2 \tau)) \Delta \bar{B}(s-\tau)\Big|dB(s)\nonumber \\
	&=\left|\xi(0)-D(\xi(-\tau))\right|^{p}+\sum_{i=1}^{7}H_{i}(t).
\end{align}
By assumption \textbf{(A5)}, \eqref{yuan B2} and  Young's inequality $a^{p-2}b\le\frac{p-2}{p}a^p+\frac{2}{p}b^{p-2}$, it derives that
\begin{align}\label{H1+H2}
	\mathbb{E}&\big(\sup_{0\leq t\leq T}H_{1}(t)\big)+\mathbb{E}\big(\sup_{0\leq t\leq T}H_{2}(t)\big)\nonumber\\
	=&\mathbb{E}\Big(\sup_{0\leq t\leq T}\int_{0}^{t}p|y(s)-D(\bar{y}(s-\tau))|^{p-2}\Big[(\bar{y}(s)-D(\bar{y}(s-\tau)))^Tb_h(\bar{y}(s),\bar{y}(s-\tau))\nonumber\\
	 &+(p-1)|\sigma(\bar{y}(s),\bar{y}(s-\tau))|^2\Big]ds\Big)\nonumber\\
	\leq &\mathbb{E}\int_{0}^{T}p|y(s)-D(\bar{y}(s-\tau))|^{p-2}\Big[(\bar{y}(s)-D(\bar{y}(s-\tau)))^Tb_h(\bar{y}(s),\bar{y}(s-\tau))\nonumber\\
	&+(p-1)|\sigma(\bar{y}(s),\bar{y}(s-\tau))|^2\Big]ds	\nonumber\\
	\leq &C\mathbb{E}\int_{0}^{T}|y(s)-D(\bar{y}(s-\tau))|^{p-2}(1+|\bar{y}(s)|^2+|\bar{y}(s-\tau)|^2)ds\nonumber\\
	\leq &C\mathbb{E}\int_{0}^{T}\big[|y(s)|^{p-2}+|D(\bar{y}(s-\tau))|^{p-2}\big](1+|\bar{y}(s)|^2+|\bar{y}(s-\tau)|^2)ds\nonumber\\
	\leq &C\mathbb{E}\int_{0}^{T}\big[1+|y(s)|^{p}+|\bar{y}(s)|^p+|\bar{y}(s-\tau)|^p\big]ds\nonumber \\
	\leq &C+C\int_{0}^{T}\mathbb{E}\sup _{0\leq u\leq s}\big|y(u)\big|^pds.
\end{align}
By Young's inequality and assumption \textbf{(A4)}, it derives that
\begin{align}\label{H3}
	\mathbb{E}&\big(\sup_{0\leq t\leq T}H_3(t)\big)\nonumber\\
	=&\mathbb{E}\Big(\sup_{0\leq t\leq T} \int_{0}^{t}p(p-1)|y(s)-D(\bar{y}(s-\tau))|^{p-2}|\sigma_1(\bar{y}(s),\bar{y}(s-\tau))\sigma(\bar{y}(s),\bar{y}(s-\tau))\Delta\bar{B}(s)\nonumber\\
	&+\sigma_2(\bar{y}(s),\bar{y}(s-\tau))\sigma(\bar{y}(s-\tau),\bar{y}(s-2\tau))\Delta\bar{B}(s-\tau)|^2ds\Big)\nonumber\\
	\leq &  C\mathbb{E}\int_{0}^T|y(s)-D(\bar{y}(s-\tau))|^pds+C\mathbb{E}\int_{0}^{T}|\sigma_1(\bar{y}(s),\bar{y}(s-\tau))\sigma(\bar{y}(s),\bar{y}(s-\tau))\Delta\bar{B}(s)\nonumber\\
	&+\sigma_2(\bar{y}(s),\bar{y}(s-\tau))\sigma(\bar{y}(s-\tau),\bar{y}(s-2\tau))\Delta\bar{B}(s-\tau)|^pds\nonumber\\
	\leq &  C\mathbb{E}\int_{0}^T|y(s)-D(\bar{y}(s-\tau))|^pds+C\mathbb{E}\int_{0}^{T}|\sigma_1(\bar{y}(s),\bar{y}(s-\!\tau))\sigma(\bar{y}(s),\bar{y}(s-\!\tau))\Delta\bar{B}(s)|^pds\nonumber\\
	+&C\mathbb{E}\int_{0}^{T}|\sigma_2(\bar{y}(s),\bar{y}(s-\tau))\sigma(\bar{y}(s-\tau),\bar{y}(s-2\tau))\Delta\bar{B}(s-\tau)|^pds\nonumber\\
	\leq &  C\mathbb{E}\int_{0}^T|y(s)-D(\bar{y}(s-\tau))|^pds+C\mathbb{E}\int_{0}^{T}|\sigma_1(\bar{y}(s),\bar{y}(s\!-\!\tau))\sigma(\bar{y}(s),\bar{y}(s\!-\!\tau))|^pds\nonumber\\
	+&C\mathbb{E}\!\int_{0}^{T}|\Delta\bar{B}(s)|^pds\!+\!C\mathbb{E}\!\int_{0}^{T}\big(|\sigma_2(\bar{y}(s),\bar{y}(s\!-\!\tau))\sigma(\bar{y}(s\!-\!\tau),\bar{y}(s\!-\!2\tau))|^p\!+\!|\Delta\bar{B}(s\!-\!\tau)|^p\big)ds\nonumber\\
	\leq &C\mathbb{E}\int_{0}^T|y(s)-D(\bar{y}(s-\tau))|^pds+C\mathbb{E}\int_{0}^{T}\big|(|(\bar{y}(s)|+|\bar{y}(s-\tau)|)\big|^pds+T\Delta^{\frac{p}{2}}\nonumber\\
	+&C\mathbb{E}\int_{0}^{T}\big|(|(\bar{y}(s-\tau)|+|\bar{y}(s-2\tau)|)\big|^pds+T\Delta^{\frac{p}{2}}\nonumber\\
	\leq &C_p+C_p\mathbb{E}(\sup_{0\leq t\leq T}|y(t)-D(\bar{y}(t-\tau))|^{p})+C_p\int_{0}^{T}\mathbb{E}(\sup_{0\leq u\leq s}|y(u)|^p)ds.
	\end{align}
By Young's inequality and \eqref{bh} then
\begin{align}\label{H4}
	\mathbb{E}&\big(\sup_{0\leq t\leq T}H_{4}(t)\big)\nonumber\\
	=&p\mathbb{E}\Big(\sup_{0\leq t\leq T}\int_{0}^{t}|y(s)-D(\bar{y}(s-\tau))|^{p-2}(y(s)-\bar{y}(s))^Tb_{h}(\bar{y}(s),\bar{y}(s-\tau))ds
\Big)\nonumber\\
	\leq & p\mathbb{E}\int_{0}^{T}|y(s)-D(\bar{y}(s-\tau))|^{p-2}(y(s)-\bar{y}(s))^Tb_{h}(\bar{y}(s),\bar{y}(s-\tau))ds
	\nonumber \\
	\leq &p\mathbb{E}(\sup_{0\leq t\leq T}|y(t)-D(\bar{y}(t-\tau))|^{p-2}[\mathbb{E}(\int_{0}^{T}|y(s)-\bar{y}(s)|| b_{h}(\bar{y}(s),\bar{y}(s-\tau))|ds)^{1/2}]^2)\nonumber \\
	\leq& (p\!-\!2)\mathbb{E}(\sup_{0\leq t\leq T}|y(t)\!-\!D(\bar{y}(t-\tau))|^{p})+2\mathbb{E}\int_{0}^{T}|y(s)-\bar{y}(s)|^{p/2}| b_{h}(\bar{y}(s),\bar{y}(s-\tau))|^{p/2}ds\nonumber\\
	\leq& (p-2)\mathbb{E}(\sup_{0\leq t\leq T}|y(t)-D(\bar{y}(t-\tau))|^{p})+2\mathbb{E}\int_{0}^{T}|y(s)-\bar{y}(s)|^{p}\Delta^{-\alpha p}ds\nonumber\\
	\leq& (p-2)\mathbb{E}(\sup_{0\leq t\leq T}|y(t)-D(\bar{y}(t-\tau))|^{p})+2\int_{0}^{T}\mathbb{E}(\sup_{0\leq u\leq s}|y(u)|^p\Delta^{-\alpha p})ds.\nonumber\\
			\leq & C_{p}\mathbb{E}(\sup_{0\leq t\leq T}|y(t)-D(\bar{y}(t-\tau))|^{p})+C_p
		\int_{0}^{T}\mathbb{E}(\sup_{0\leq u\leq s}|y(u)|^p)ds.
\end{align}
By the BDG inequality and Young's inequality, it derives that
\begin{align}\label{H5}
	\mathbb{E}&(\sup_{0\leq t\leq T}H_{5}(t))\nonumber\\
	=& \mathbb{E}\Big(\sup_{0\leq t\leq T}p\int_{0}^{t}|y(s)-D(\bar{y}(s-\tau))|^{p-2}\langle y(s)-D(\bar{y}(s-\tau)),\sigma(\bar{y}(s),\bar{y}(s-\tau))dB(s)\rangle \Big)\nonumber\\
	\leq &C_p\mathbb{E}\Big(\int_{0}^{T}(|y(s)-D(\bar{y}(s-\tau))|^{2p-2}|\sigma(\bar{y}(s),\bar{y}(s-\tau))|^2)ds\Big)^{1/2}\nonumber\\
	\leq &C_p\mathbb{E}\Big(\sup_{0\leq t \leq T}|y(t)-D(\bar{y}(t-\tau))|^{p-1}(\int_{0}^{T}|\sigma(\bar{y}(s),\bar{y}(s-\tau))|^2ds)^{1/2}\Big)\nonumber\\
	\leq & C_p+C_p\mathbb{E}(\sup_{0\leq t\leq T}|y(t)-D(\bar{y}(t-\tau))|^{p})+C_p\int_{0}^{T}\mathbb{E}(\sup_{0\leq u\leq s}|y(u)|^p)ds.
\end{align}
Using the similar approach as in \eqref{H5}, we can estimate $H_{6}(t)$ as
\begin{align}\label{H6}
	\mathbb{E}&\big(\sup_{0\leq t\leq T}H_{6}(t)\big)\nonumber\\
	=&\mathbb{E}\Big(\sup_{0\leq t\leq T}p\int_{0}^{t}|y(s)-D(\bar{y}(s-\tau))|^{p-1}\sigma_1(\bar{y}(s),\bar{y}(s-\tau))\sigma(\bar{y}(s),\bar{y}(s-\tau))\Delta\bar{B}(s)\Big)\nonumber\\
	\leq &C_p\mathbb{E}\Big(\!\sup_{0\leq t\leq T}\!|y(t)\!-\!D(\bar{y}(t-\!\tau))|^{p-1}(\int_{0}^{T}\!|\sigma_{1}(\bar{y}(s), \bar{y}(s-\tau)) \sigma(\bar{y}(s), \bar{y}(s-\tau)) \Delta \bar{B}(s)|^2ds\big)^{1/2}\Big)\nonumber \\
	\leq & C_p+C_p\mathbb{E}(\sup_{0\leq t\leq T}|y(t)-D(\bar{y}(t-\tau))|^{p})+C_p\int_{0}^{T}\mathbb{E}(\sup_{0\leq u\leq s}|y(u)|^p)ds.
\end{align}
Using the similar approach as in \eqref{H6}, we can estimate $H_{7}(t)$ as
\begin{align}\label{H7}
	\mathbb{E}\big(\sup_{0\leq t\leq T}\!H_{7}(t)\big)\leq  C_p\!+C_pE(\sup_{0\leq t\leq T}|y(t)\!-D(\bar{y}(t\!-\!\tau))|^{p})\!+C_p\int_{0}^{T}\mathbb{E}(\sup_{0\leq u\leq s}\!|y(u)|^p)ds.
\end{align}
Now, substituting \eqref{H1+H2} \eqref{H3}\eqref{H4} \eqref{H5} \eqref{H6} and \eqref{H7} into \eqref{sum of H}, it derives that
\begin{equation}
		\mathbb{E}(\sup _{0\leq t \leq T}\big|y(t)-D(\bar{y}(t-\tau))\big|^p)\leq C+C\int_{0}^{T}\mathbb{E}\Big(\sup _{0\leq u\leq s}\big|y(u)\big|^p\Big)ds.
	\end{equation}
Then, substituting into \eqref{yt}
	\begin{equation*}
		\mathbb{E}(\sup _{0\leq t\leq T}|y(t)|^{p})\leq C+C\int_{0}^{T}\mathbb{E}\Big(\sup _{0\leq u\leq s}\big|y(u)\big|^p\Big)ds.
	\end{equation*}
An application of Gronwall inequality yields 
\begin{equation*}
	\sup_{0\le t\le T}\mathbb{E}\left|y\left(t\right)\right|^p\le C.
\end{equation*} 
\end{Proof}

\begin{lemma}
	Under \eqref{yuan B2}, then $\forall t>0$, $p \geq 2,$
	\begin{align}\label{y(t)-bary(t)}
\mathbb{E}\big[\sup_{t_k\leq t\leq t_{k+1}}|y(t)-\bar{y}(t)|^{p}\big] \leq C,
\end{align}
where the positive constant $C:=C(p,T,K_R,||\xi||).$
\end{lemma}
\begin{Proof}
\rm By $\left|\sum_{i=1}^{n} x_{i}\right|^{p} \leq n^{p-1} \sum_{i=1}^{n}\left|x_{i}\right|^{p}$ it derives that
\begin{align}\label{sum of J}
	|y(t)-\bar{y}(t)|^{p}&=\left|y(t)-y\left(t_{k}\right)\right|^{p}
	\nonumber\\
	&\leq 
	5^{p-1}\left[\left|D(y(t-\tau))-D\left(y\left(t_{k}-\tau\right)\right)\right|^{p}\nonumber
	+\left|\int_{t_{k}}^{t} b_{h}(\bar{y}(s), \bar{y}(s-\tau)) ds\right|^{p}\right.
	\nonumber\\
	&+\left|\int_{t_{k}}^{t} \sigma(\bar{y}(s), \bar{y}(s-\tau)) d B(s)\right|^{p}
	\nonumber\\
	&+\left|\int_{t_{k}}^{t} \sigma_{1}(\bar{y}(s), \bar{y}(s-\tau)) \sigma(\bar{y}(s), \bar{y}(s-\tau)) \Delta \bar{B}(s) d B(s)\right|^{p}
	\nonumber\\
	&\left.+\left|\int_{t_{k}}^{t} \sigma_{2}(\bar{y}(s), \bar{y}(s-\tau)) \sigma(\bar{y}(s-\tau), \bar{y}(s-2 \tau)) \Delta \bar{B}(s-\tau) d B(s)\right|^{p}\right]\nonumber\\
	&=5^{p-1}\Big[\left|D(y(t-\tau))-D\left(y\left(t_{k}-\tau\right)\right)\right|^{p}+J_{1}(t)+J_{2}(t)+J_{3}(t)+J_{4}(t)\Big].
\end{align}
By the H$\ddot{\boldsymbol{\text{o}}}$lder inequality and \eqref{bh}, it derives that
\begin{align}\label{J1}
	\mathbb{E}\Big[\sup_{t_k\leq t\leq t_{k+1}}(J_{1}(t))\Big]=&\mathbb{E}\Big[\sup_{t_k\leq t\leq t_{k+1}}\left|\int_{t_{k}}^{t} b_{h}(\bar{y}(s), \bar{y}(s-\tau)) ds\right|^{p}\Big]\nonumber\\
	\leq & \Delta^{p-1} \mathbb{E}\Big[\int_{t_{k}}^{t_{k+1}} \mid b_{h}\left(\bar{y}(s),\left.\bar{y}(s-\tau)\right|^{p} d s\right.\Big]\nonumber\\
	\leq & K_R\Delta^{p-1} \mathbb{E} \Big[\int_{t_{k}}^{t_{k+1}} |\bar{y}(s)|^p+|\bar{y}(s-\tau)|^{p} d s\Big]\nonumber\\
	\leq &K_R\Delta^{p-1}C_p\Big[\tau||\xi||^p+\mathbb{E}\int_{t_k}^{t}\sup_{0\leq u\leq s}|y(u)|^{p}ds\Big]\nonumber\\
	\leq & C+C\mathbb{E}\int_{t_k}^{t}\sup_{0\leq u\leq s}|y(u)|^{p}ds.
\end{align}
By the BDG inequality, H$\ddot{\boldsymbol{\text{o}}}$lder inequality and \textbf{(A3)}, it derives that
\begin{align}\label{J2}
	\mathbb{E}\Big[\sup_{t_k\leq t\leq t_k+1}(J_{2}(t))\Big]=&\mathbb{E}\Big[\sup_{t_k\leq t\leq t_{k+1}}\left|\int_{t_{k}}^{t} \sigma(\bar{y}(s), \bar{y}(s-\tau)) d B(s)\right|^{p}\Big]\nonumber\\
	\leq & C_p\mathbb{E} \Big[\int_{t_{k}}^{t_{k+1}}|\sigma(\bar{y}(s), \bar{y}(s-\tau))|^{2} d s\Big]^{\frac{p}{2}}\nonumber\\
	\leq &C_p\Delta^{\frac{p-2}{2}}\mathbb{E}\Big[\int_{t_{k}}^{t_{k+1}}|\sigma(\bar{y}(s), \bar{y}(s-\tau))|^{p} d s\Big]\nonumber\\
	\leq &K_RC_p\Delta^{\frac{p-2}{2}}\mathbb{E}\Big[\int_{t_{k}}^{t_{k+1}}\big(1+|\bar{y}(s)|^p+|\bar{y}(s-\tau)|^p\big) d s\Big]\nonumber\\
	\leq &K_R\Delta^{\frac{p-2}{2}}C_p\Big[\tau||\xi||^p+\mathbb{E}\int_{t_k}^{t}\sup_{0\leq u\leq s}|y(u)|^pds\Big]\nonumber\\
	\leq &C+C\mathbb{E}\int_{t_k}^{t}\sup_{0\leq u\leq s}|y(u)|^pds.
\end{align}
Using the similar approach as in \eqref{J2}, it estimates $J_{3}(t)$ and $J_{4}(t)$ as
\begin{align}\label{J3 AND J4}
	\mathbb{E}&\Big[\sup_{t_k\leq t\leq t_k+1}(J_{3}(t))\Big]+\mathbb{E}\Big[\sup_{t_k\leq t\leq t_k+1}(J_{4}(t))\Big]\nonumber\\
	=&\mathbb{E}\Bigg[\sup_{t_k\leq t\leq t_{k+1}}\left|\int_{t_{k}}^{t} \sigma_{1}(\bar{y}(s), \bar{y}(s-\tau)) \sigma(\bar{y}(s), \bar{y}(s-\tau)) \Delta \bar{B}(s) d B(s)\right|^{p}\Bigg] \nonumber\\
	+&\mathbb{E}\Bigg[\sup_{t_k\leq t\leq t_{k+1}}\left|\int_{t_{k}}^{t} \sigma_{2}(\bar{y}(s), \bar{y}(s-\tau)) \sigma(\bar{y}(s-\tau), \bar{y}(s-2\tau)) \Delta \bar{B}(s) d B(s-\tau)\right|^{p}\Bigg] \nonumber\\
	\leq &C_{p} \Delta^{\frac{p-2}{2}} \mathbb{E} \int_{t_{k}}^{t_{k+1}}|\sigma_{1}(\bar{y}(s), \bar{y}(s-\tau)) \sigma(\bar{y}(s), \bar{y}(s-\tau)) \Delta \bar{B}(s)|^{p} d s\nonumber\\
	+&C_{p}\Delta^{\frac{p-2}{2}} \mathbb{E} \int_{t_{k}}^{t_{k+1}}\left|\sigma_{2}(\bar{y}(s), \bar{y}(s-\tau)) \sigma(\bar{y}(s-\tau), \bar{y}(s-2 \tau)) \Delta \bar{B}(s-\tau)\right|^{p} d s\nonumber\\
	\leq & K_RC_p\Delta^{\frac{p-2}{2}}\Big[\tau||\xi||^p+\int_{t_k}^{t}\mathbb{E} \sup_{0\leq u\leq s}|y(u)|^pds+(\Delta\bar{B}(s))^pds\Big]\nonumber\\
	\leq &C+C\int_{t_k}^{t}\mathbb{E} \sup_{0\leq u\leq s}|y(u)|^pds.
\end{align}
Now, substituting \eqref{J1} \eqref{J2} and \eqref{J3 AND J4} into \eqref{sum of J}, it derives that
\begin{align*}
	\mathbb{E}|y(t)-\bar{y}(t)|^{p} \leq C+\int_{t_k}^{t}\mathbb{E} \sup_{0\leq u\leq s}|y(u)|^pds.
\end{align*}
By \eqref{y(t)}, we could obtain 
\begin{align*}
	\mathbb{E}\big[\sup_{t_k\leq t\leq t_{k+1}}|y(t)-\bar{y}(t)|^{p}\big] \leq C.
\end{align*}
\end{Proof}

\begin{lemma}
	If assumptions \rm\textbf{(A2)} \emph and \textbf{(A5)} \emph{hold, for any real number $R>\|\xi\|$ define a stopping time $\tau_{R}=\inf \{t \geq 0:|x(t)| \geq R\}$ and let $\inf \ \Phi=\infty$ then
	\begin{equation}\label{P1}
P\left(\tau_{R} \leq T\right) \leq \frac{C}{R^{2}}.
\end{equation}
where the positive constant $C:=C(T,K_1,||\xi||,\kappa).$}
\end{lemma}
\begin{Proof}
\rm By the same way as in \eqref{basic inequality} and $
\text { Itô }
$ formula, it derives that
\begin{align}\label{stoping time 1}
	\mathbb{E}\left|x\left(t \wedge \tau_{R}\right)\right|^{2}\leq &\frac{1}{1-\kappa} \mathbb{E}\left|x\left(t \wedge \tau_{R}\right)-D\left(x\left(t \wedge \tau_{R}-\tau\right)\right)\right|^{2}+\frac{1}{\kappa} \mathbb{E}\left|D\left(x\left(t \wedge \tau_{R}-\tau\right)\right)\right|^{2} \nonumber \\
 \leq &\frac{1}{1\!-\!\kappa}\Big[\big|\xi(0)\!-\!D(\xi(\!-\!\tau))\big|^{2}\!+\!\mathbb{E} \int_{0}^{t \wedge \tau_{R}} 2x(s)\!-\!D(x(s\!-\!\tau)))^{T} \cdot (b(x(s), x(s\!-\!\tau))\nonumber \\
&+\big|\sigma(x(s), x(s-\tau))\big|^{2} d s\Big]+\frac{1}{\kappa} \mathbb{E}\left|D\left(x\left(t \wedge \tau_{R}-\tau\right)\right)\right|^{2}\nonumber \\
=&\frac{1}{1-\kappa}\Big[\big|\xi(0)-D(\xi(-\tau))\big|^2+ L_{1}+L_{2}\Big]+L_{3}.
\end{align}
By Assumption \textbf{(A5)}, it derives that
\begin{align}\label{L1+L2}
	L_{1}+L_{2}=&\mathbb{E} \int_{0}^{t \wedge \tau_{R}} 2(x(s)-D(x(s-\tau)))^{T} \cdot b(x(s), x(s-\tau))+\big|\sigma(x(s), x(s-\tau))\big|^{2} d s\Big]\nonumber\\
	\leq &\mathbb{E}\int_{0}^{t\wedge \tau_{R}}2K_{1}\big(1+|x(s)|^2+|x(s-\tau)|\big)ds\nonumber\\
	\leq &2K_{1} T+2K_{1}\tau||\xi||^2+4K_{1} \int_{0}^{t} \mathbb{E}|x(s\wedge \tau_{R})|^2ds.
\end{align}
By Assumption \textbf{(A2)}, it derives that
\begin{align}\label{L3}
L_{3}=&\frac{1}{\kappa}\mathbb{E}|D(x(t\wedge \tau_{R}-\tau))|^2\nonumber\\
\leq &\kappa\mathbb{E}|x(t\wedge \tau_{R}-\tau)|^2\nonumber\\
\leq &\kappa\Big[||\xi||^2+E(x(t\wedge\tau_{R}))^2\Big].
\end{align}
Now, substituting \eqref{L1+L2} and \eqref{L3} into \eqref{stoping time 1}, it derives that
\begin{align*}
\mathbb{E}|x(t\wedge \tau_{R})|^2\leq C+C\int	_{0}^{t} \mathbb{E}|x(s\wedge \tau_{R})|^2ds.
\end{align*}
The Gronwall inequality yields that
\begin{align*}
	\mathbb{E}|x(t\wedge \tau_{R})|^2\leq C.
\end{align*}
This implies
\begin{align*}
	R^2\cdot P(\tau_{R}\leq T)\leq C.
\end{align*}
\end{Proof}
\begin{lemma}
	If assumptions \rm\textbf{(A1)-(A3)} \emph{hold, for any real number $R>\|\xi\|$ define a stopping time $\rho_{R}=\inf \{t \geq 0:|x(t)| \geq R\}$ and let $\inf\ \Phi=\infty$ then
	\begin{equation}\label{P2}
P\left(\rho_{R} \leq T\right) \leq \frac{C}{R^{2}}.
\end{equation}}
\end{lemma}
\begin{Proof}
\rm Using the similar approach as in \eqref{basic inequality}, we could obtain
\begin{align}\label{stopping time 2}
	\mathbb{E}\left|y\left(t \wedge \rho_{R}\right)\right|^{2} \leq &\frac{1}{1-\kappa} \mathbb{E}\left|y\left(t \wedge \rho_{R}\right)-D\left(\bar{y}\left(t \wedge \rho_{R}-\tau\right)\right)\right|^{2}+\frac{1}{\kappa} \mathbb{E}\left|D\left(\bar{y}\left(t \wedge \rho_{R}-\tau\right)\right)\right|^{2} \nonumber \\
	\leq &\frac{1}{1-\kappa}\Big[|\xi(0)-D(\xi(-\tau))|^{2}\nonumber \\
	+&2 \mathbb{E}\int_{0}^{t \wedge \rho_{R}}\Big((\bar{y}(s)\!-\!D(\bar{y}(s\!-\!\tau)))^{T} \cdot b_{h}(\bar{y}(s), \bar{y}(s\!-\!\tau)+|\sigma(\bar{y}(s), \bar{y}(s-\tau))|^{2}\Big) d s\nonumber \\
	+&2 \mathbb{E} \int_{0}^{t \wedge \rho_{R}} \mid \sigma_{1}(\bar{y}(s), \bar{y}(s-\tau)) \cdot \sigma(\bar{y}(s), \bar{y}(s-\tau)) \Delta \bar{B}(s)\nonumber \\
	&\quad+\left.\sigma(\bar{y}(s-\tau), \bar{y}(s-2 \tau)) \Delta \bar{B}(s-\tau)\right|^{2} d s\nonumber \\
	&+2 \mathbb{E} \int_{0}^{t \wedge \rho_{R}}(y(s)-\bar{y}(s))^{T} \cdot b_{h}(\bar{y}(s), \bar{y}(s-\tau) ds\Big]\nonumber \\
	&+\frac{1}{\kappa} \mathbb{E}\left|D\left(\bar{y}\left(t \wedge \rho_{R}-\tau\right)\right)\right|^{2}\nonumber \\
	=&\frac{1}{1-\kappa}\Big[|\xi(0)-D(\xi(-\tau))|^2+M_{1}+M_{2}+M_{3}+M_{4}\Big]+M_{5}.
\end{align}
By \eqref{yuan B2} and \eqref{y(t)-bary(t)} it derives that
\begin{align}\label{M1+M2}
	M_{1}+M_{2}&=2 \mathbb{E}\int_{0}^{t \wedge \rho_{R}}\Big((\bar{y}(s)-D(\bar{y}(s-\tau)))^{T} \cdot b_{h}(\bar{y}(s), \bar{y}(s-\tau)+|\sigma(\bar{y}(s), \bar{y}(s-\tau))|^{2}\Big) d s\nonumber \\
	&\leq 2K_{2}\mathbb{E} \int_{0}^{T }(1+\sup_{0\leq u\leq s}|\bar{y}(u\wedge \rho_{R})|^2+\sup_{0\leq u\leq s}|\bar{y}(u\wedge \rho_{R}-\tau|^2)ds\nonumber \\
	&\leq2TK_{2}+2K_{2}\tau ||\xi||^2+4K_{2}\mathbb{E}\int_{0}^{T }\sup_{0\leq u\leq s}|y(u\wedge\rho_{R})|^2du\nonumber \\
	&\leq C+C\mathbb{E}\int_{0}^{T}\sup_{0\leq u\leq s}|y(u\wedge\rho_{R})|^2du.
\end{align}
Using the similar approach as in \eqref{J3 AND J4}, we can estimate $M_{3}(t)$ as
\begin{align}\label{M3}
	M_{3}=&2 \mathbb{E} \int_{0}^{t \wedge \rho_{R}} \mid \sigma_{1}(\bar{y}(s), \bar{y}(s-\tau)) \cdot \sigma(\bar{y}(s), \bar{y}(s-\tau)) \Delta \bar{B}(s)\nonumber \\
	&+\sigma_2(\bar{y}(s), \bar{y}(s-\tau))\left.\sigma(\bar{y}(s-\tau), \bar{y}(s-2 \tau)) \Delta \bar{B}(s-\tau)\right|^{2} d s\nonumber \\
	\leq &C+C\int_{0}^{T}\mathbb{E} \sup_{0\leq u\leq s}|y(u\wedge\rho_{R})|^2ds.
\end{align}
By \eqref{bh1} and \eqref{y(t)-bary(t)} it derives that
\begin{align}\label{M4}
	M_{4}&=2 \mathbb{E} \int_{0}^{t \wedge \rho_{R}}(y(s)-\bar{y}(s))^{T} \cdot b_{h}(\bar{y}(s), \bar{y}(s-\tau) ds\nonumber\\
	&\leq 2K_2C\mathbb{E}\int_{0}^{T}.(1+\sup_{0\leq u\leq s}|\bar{y}(u\wedge \rho_{R})|^2+\sup_{0\leq u\leq s}|\bar{y}(u\wedge \rho_{R}-\tau|^2)ds\nonumber \\
	&\leq2TK_{2}+2K_{2}\tau ||\xi||^2+4K_{2}\mathbb{E}\int_{0}^{T }\sup_{0\leq u\leq s}|y(u\wedge\rho_{R})|^2du\nonumber \\
	&\leq C+C\mathbb{E}\int_{0}^{T}\sup_{0\leq u\leq s}|y(u\wedge\rho_{R})|^2du.
\end{align}
By assumption \textbf{(A3)}, it derives that
\begin{align}\label{M5}
	M_{5}&=\frac{1}{\kappa}\mathbb{E}|D(\bar{y}\wedge\rho_R-\tau)|\nonumber\\
	&\leq \kappa\mathbb{E}|\bar{y}(t\wedge\rho_{R}-\tau)|^{2}\nonumber\\
	&\leq \kappa \sup_{0\leq u \leq t}\mathbb{E}|\bar{y}(u\wedge\rho_{R}-\tau)|^{2}\nonumber\\
	&\leq \kappa(||\xi||^2+\sup_{0\leq u\leq t}\mathbb{E}|\bar{y}(u\wedge\rho_{R})|^{2})\nonumber\\
	&\leq C+C\sup_{0\leq u\leq t}\mathbb{E}|y(u\wedge\rho_{R})|^{2}).
\end{align}
Now, substituting \eqref{M1+M2} \eqref{M3} \eqref{M4} and \eqref{M5}  into \eqref{stopping time 2}
\begin{align*}
\mathbb{E}|y(t\wedge \rho_{R})|^2\leq C+C\int	_{0}^{t} \mathbb{E}|y(u\wedge \rho_{R})|^2du.
\end{align*}
An application of Gronwall inequality yeilds:
\begin{align*}
	\mathbb{E}|y(t\wedge \rho_{R})|^2\leq C.
\end{align*}
This implies
\begin{align*}
	R^2\cdot P(\rho_{R}\leq T)\leq C.
\end{align*}
\end{Proof}
\subsection{Proof of Therom 4.1}
In this section, we give proof of the main theorem of this paper, the strong convergence of the tamed Milstein\eqref{tamed method} to the solution of \eqref{NSDDEs1}.
\begin{Proof}
	\rm Let $\tau_R$ and $ \rho_R$ be the same as before, define: $\theta_R=\tau_R\wedge \rho_R$ and $e\left(t\right)=x\left(t\right)-y\left(t\right)$, by Young's inequality, it derives that
	\begin{equation}
\begin{aligned}
\mathbb{E}\left[\sup _{0 \leq t \leq T}|e(t)|^{p}\right]\leq &\mathbb{E}\left[\sup _{0 \leq t \leq T}|e(t)|^{p} I_{\left\{\tau_{R} \leq T \text { or } \rho_{R} \leq T\right\}}\right]+\mathbb{E}\left[\sup _{0 \leq t \leq T}\left|e\left(t \wedge \theta_{R}\right)\right|^{p}\right] \\
\leq &\frac{p \eta}{q} \mathbb{E}\left[\sup _{0 \leq t \leq T}|x(t)-y(t)|^{q}\right]+\frac{q-p}{q \eta^{p /(q-p)}} P\left(\tau_{R} \leq T \text { or } \rho_{R} \leq T\right)\nonumber\\
&+\mathbb{E}\left[\sup _{0 \leq t \leq T}\left|e\left(t \wedge \theta_{R}\right)\right|^{p}\right].
\end{aligned}
\end{equation}
Using the similar approach as \eqref{basic inequality}, it derives that 
\begin{align}\label{et}
	&\mathbb E\big(\sup_{0\leq t\leq T}|e(t\wedge\theta_R)|^p\big)\nonumber\\
	&\leq\frac{1}{(1\!-\!\kappa)^p}\mathbb E\big(\sup_{0\leq t\leq T}\!|\!x(t\wedge\theta_R)\!-\!y(t\wedge\theta_R)\!-\!D(x(t\wedge\theta_R-\tau))\!+\!D(\bar{y}(t\wedge\theta_R\!-\!\tau))|^p\big)\!+\!C.
\end{align}
By the H$\ddot{\boldsymbol{\text{o}}}$lder inequality and \textbf{(A3)}, for any $t \in [0,T]$, it derives that
\begin{align}\label{sum of Q}
	&|x(t\wedge\theta_R)-y(t\wedge\theta_R)-D(x(t\wedge\theta_R-\tau))+D(\bar{y}(t\wedge\theta_R-\tau))|^p\nonumber\\
	&\leq C_R\int_{0}^{t\wedge\theta_R}\Big[|(x(s)-y(s)|^p+|x(s-\tau)-y(s-\tau)|^p\Big]ds\nonumber\\
	&+C\int_{0}^{t\wedge\theta_R}|b(y(s),y(s-\tau))-b(\bar{y}(s),\bar{y}(s-\tau))|^pds\nonumber\\
	&+C\int_{0}^{t\wedge\theta_R}|b(\bar{y}(s),\bar{y}(s-\tau))-b_h(\bar{y}(s),\bar{y}(s-\tau))|^pds\nonumber\\
	&+C|\int_{0}^{t\wedge\theta_R}(\sigma(x(s),x(s-\tau))-\sigma(\bar{y}(s),\bar{y}(s-\tau))dB(s)|^p\nonumber\\
	&+C|\int_{0}^{t\wedge\theta_R}(\sigma(\bar{y}(s),\bar{y}(s-\tau))\sigma_1(\bar{y}(s),\bar{y}(s-\tau))dB(s)|^p\nonumber\\
	&+C|\int_{0}^{t\wedge\theta_R}\sigma_2(\bar{y}(s),\bar{y}(s-\tau))\sigma(\bar{y}(s),\bar{y}(s-2\tau))dB(s))|^p\nonumber\\
	&=\sum_{i=1}^{6}Q_i(t).
\end{align}
We now estimate $Q_1(t)$, by \textbf{(A1)} it derives that
\begin{align}\label{Q1}
	\mathbb{E}&(\sup_{0\leq t\leq T}Q_{1}(t))\nonumber\\
	=&\mathbb{E}(\sup_{0\leq t\leq T} C_R\int_{0}^{t\wedge\theta_R}\Big[|(x(s)-y(s)|^p+|x(s-\tau)-y(s-\tau)|^p\Big]ds)\nonumber\\
	\leq &\mathbb{E}\int_{0}^{T}\sup_{0\leq u\leq s}\Big[|x(u\wedge\theta_R)-y(u\wedge\theta_R)|^p+|x(u\wedge\theta_R-\tau))-y(u\wedge\theta_R-\tau))|^p\Big]ds\nonumber\\
	\leq &\mathbb{E}\int_{0}^{T}\sup_{0\leq u\leq s}|e(u\wedge\theta_R)|^pds+C\Delta^p.
	\end{align}
By \textbf{(A3)} and \eqref{y(t)-bary(t)}, it derives that
\begin{align}\label{Q2}
	\mathbb{E}&(\sup_{0\leq t \leq T}Q_2(t))\nonumber\\
	=&\mathbb{E}(\sup_{0\leq t \leq T}C\int_{0}^{t\wedge\theta_R}|b(y(s),y(s-\tau))-b(\bar{y}(s),\bar{y}(s-\tau))|^pds)\nonumber\\
	\leq &\mathbb{E}\int_{0}^{T}|y(s\wedge\theta_R)-\bar{y}(s\wedge\theta_R)|^pds+C\int_{0}^{T}|\xi(s\wedge\theta_R)-\xi(\kappa(s\wedge\theta_R))|^pds\nonumber\\
	\leq & C\Delta^{\frac{p}{2}}.
\end{align}
According to \textbf{(B1)}, we could obtain
\begin{align}\label{Q3}
	\mathbb{E}&(\sup_{0\leq t \leq T}Q_3(t))\nonumber\\
	=&\mathbb{E}(\sup_{0\leq t \leq T}C\int_{0}^{t\wedge\theta_R}|b(\bar{y}(s),\bar{y}(s-\tau))-b_h(\bar{y}(s),\bar{y}(s-\tau))|^pds)\nonumber\\
	\leq &C\mathbb{E}\int_{0}^{T}\Big|b(\bar{y}(s\wedge\theta_R),\bar{y}(s\wedge\theta_R-\tau))-b_h(\bar{y}(s\wedge\theta_R),\bar{y}(s\wedge\theta_R-\tau))\Big|^pds\nonumber\\
	\leq &C{N_R}^p\Delta^p\leq C.
\end{align}
Using similar approach as in \eqref{I3} \eqref{Q1} and \eqref{Q2} we have
\begin{align}\label{Q4}
	\mathbb{E}&(\sup_{0\leq t \leq T}Q_4(t))\nonumber\\
	=&\mathbb{E}(\sup_{0\leq t \leq T}C|\int_{0}^{t\wedge\theta_R}(\sigma(x(s),x(s-\tau))-\sigma(\bar{y}(s),\bar{y}(s-\tau))dB(s)|^p\nonumber\\
	\leq &C_p\Delta^{\frac{p-2}{2}}\mathbb{E}(\int_{0}^{T}|\sigma(x(s\wedge\theta_R),x(s\wedge\theta_R-\tau))-\sigma(\bar{y}(s\wedge\theta_R),\bar{y}(s\wedge\theta_R-\tau))|^pds\nonumber\\
	\leq &C_p\Delta^{\frac{p-2}{2}}\mathbb{E}\int_{0}^{T}\sup_{0\leq u\leq s}\Big[|x(u\wedge\theta_R)-y(u\wedge\theta_R)|^p+|x(u\wedge\theta_R-\tau))-y(u\wedge\theta_R-\tau))|^p\Big]ds\nonumber\\
	+&C_p\Delta^{\frac{p-2}{2}}\mathbb{E}\int_{0}^{T}|y(s\wedge\theta_R)-\bar{y}(s\wedge\theta_R)|^pds+C\int_{0}^{T}|\xi(s\wedge\theta_R)-\xi(\kappa(s\wedge\theta_R))|^pds\nonumber\\
	\leq &\mathbb{E}\int_{0}^{T}\sup_{0\leq u\leq s}|e(u\wedge\theta_R)|^pds+C\Delta^\frac{p}{2}.
	\end{align}
By \eqref{I3} and \eqref{yt} and \eqref{y(t)-bary(t)}, it derives that
\begin{align}\label{Q5}
	&\quad\mathbb{E}(\sup_{0\leq t \leq T}Q_5(t))\nonumber\\
	&\leq C\mathbb{E}\int_{0}^{T}|\sigma(\bar{y}(s\wedge\theta_R),\bar{y}(s\wedge\theta_R-\tau))\sigma_1(\bar{y}(s\wedge\theta_R),\bar{y}(s\wedge\theta_R-\tau))\Delta \bar B(s\wedge\theta_R)|^pds\nonumber\\
	&\leq C\mathbb{E}\int_{0}^{T}\Big(|\sigma(\bar{y}(s\wedge\theta_R),\bar{y}(s\wedge\theta_R-\tau))\sigma_1(\bar{y}(s\wedge\theta_R),\bar{y}(s\wedge\theta_R-\tau))\Delta \bar B(s\wedge\theta_R)\nonumber\\
	&\quad-\sigma(y(s\wedge\theta_R),y(s\wedge\theta_R-\tau))\sigma_1(y(s\wedge\theta_R),y(s\wedge\theta_R-\tau))\Delta \bar B(s\wedge\theta_R)|^p\nonumber\\
	&\quad+|\sigma(y(s\wedge\theta_R),y(s\wedge\theta_R-\tau))\sigma_1(y(s\wedge\theta_R),y(s\wedge\theta_R-\tau))\Delta \bar B(s\wedge\theta_R)|^p\Big)ds\nonumber\\
	&\leq C\mathbb{E}\int_{0}^{T}\Big(|y(s\wedge\theta_R)-\bar{y}(s\wedge\theta_R)|^p+|y(s\wedge\theta_R-\tau)-\bar{y}(s\wedge\theta_R-\tau)|^p\nonumber\\
	&\quad+|y(s\wedge\theta_R)|^p+|y(s\wedge\theta_R-\tau)|^p\Big)ds\nonumber\\
	&\quad \leq C\Delta^\frac{p}{2}+C.
\end{align}
Using the similary approach as in \eqref{Q5}, we can estimate $Q_6(t)$ as 
\begin{align}\label{Q6}
	\mathbb{E}(\sup_{0\leq t \leq T}Q_6(t))\leq C\Delta^\frac{p}{2}+C.
\end{align}
Substituting \eqref{Q1},\eqref{Q2},\eqref{Q3},\eqref{Q4},\eqref{Q5} and \eqref{Q6} into \eqref{et}
\begin{align}
		&\quad\mathbb E\big(\sup_{0\leq t\leq T}|e(t\wedge\theta_R)|^p\big)\leq C\mathbb{E}\int_{0}^{T}\sup_{0\leq u\leq s}|e(u\wedge\theta_R)|^pds+C.
\end{align}
An application of Gronwall inequality yields:
\begin{align}
	\mathbb E\big(\sup_{0\leq t\leq T}|e(t\wedge\theta_R)|^p\big)\leq C.
\end{align}
Where the positive constant $C:=C(R,||\xi||,\kappa,\Delta).$\\
By \eqref{x(t)} \eqref{y(t)} and given an $\epsilon>0$, there\ exists\ $\eta$ small enough  that
\begin{equation*}
\frac{p\eta }{q} \mathbb{E}\Big[\sup_{0\leq t\leq T}|x(t)-y(t)|^{q}\Big] \leq 2^q\frac{p\eta }{q}  \mathbb{E}\Big\{\big[\sup_{0\leq t\leq T}|x(t)|^{q}\big]+\big[\sup_{0\leq t\leq T}|y(t)|^{q}]\Big\} \leq \frac{\epsilon}{3}.
\end{equation*}
Choose R large enough that
\begin{equation*}
\frac{q-p}{q \eta^{p /(q-p)}}\left[P\left(\tau_{R} \leq T\right)+P\left(\rho_{R} \leq T\right)\right] \leq \frac{\epsilon}{3}.
\end{equation*}
And finally by Lemmas 4.1 and 4.2, we choose $\Delta$ sufficiently small, such that 
\begin{equation*}
 \mathbb E\big(\sup_{0\leq t\leq T}|e(t\wedge\theta_R)|^p\big)\leq \frac{\epsilon}{3}.
\end{equation*}
\end{Proof}

%\section*{Declarations}
%\textbf{Author Contributions:} {Methodology, Writing-Original Draft Preparation, Q.F.; Formal Analysis,Writing-Original Draft Preparation, M.Y. and J.Y.; Writing-Review and Editing, M.Y. and W.Y.; Supervision, J.Y. }\\
%\\
%\textbf{Funding:} {This research was supported by Basic Science Research Program through the National Research Foundation of Korea (NRF) funded by the Ministry of Education, Science and Technology (NRF-2016R1D1A1B03930571); National Science Foundation (DMS-1109063); ZhengJiang Social Foundation of Grant NO.20JD2D071; Zhejiang University of Science and Technology of Grant NO.F701108K03}\\
%\\
%\textbf{Institutional Review  Board Statement:} {Not applicable}\\
%\\
%\textbf{Informed Consent Statement:} {Not applicable}\\
%\\
%\textbf{Data Availability Statement:} {Data sharing not applicable to this article as no datesets were generated or analyzed during the current study}\\
%\\ 
%\textbf{Acknowledgments:} {In this section you can acknowledge any support given which is not covered by the author contribution or funding sections. This may include administrative and technical support, or donations in kind (e.g., materials used for experiments).}\\
%\\
%\textbf{Conflicts of Interest:} {The authors declare no conflict of interest} \\
%\\
%\textbf{Abbreviations:}{
%The following abbreviations are used in this manuscript:\\
%
%\noindent 
%\begin{tabular}{@{}ll}
%SDEs & Stochastic Differential Equations\\
%NSDDEs & Neutral Stochastic Delay Differential Equations\\
%\end{tabular}}

	\bibliographystyle{plain}
	
\bibliography{paper}
\end{document}